\documentclass[12pt]{article}
\usepackage{amsfonts}
\usepackage{}
\usepackage{wasysym,amssymb,eufrak,indentfirst,graphicx,cite,amsthm,color}
\usepackage[bookmarksnumbered, colorlinks, plainpages]{hyperref}
\newtheorem{Theorem}{Theorem}[section]
\newtheorem{Corollary}[Theorem]{Corollary}
\newtheorem{Lemma}[Theorem]{Lemma}
\newtheorem{Proposition}[Theorem]{Proposition}

\newtheorem{Definition}[Theorem]{Definition}
\newtheorem{Example}[Theorem]{Example}
\newtheorem{Remark}{{Remark}}

\def\s{\sigma}
\def\f{\noindent}
\def\demo{\f{ \textbf{Proof}}\hskip10pt}
\def\qed{\hfill $\Box$}
\def\mH{\mathcal{H}$$}

\begin{document}
\baselineskip 17pt

\title{New characterizations of a normal subgroup to be hypercyclically embedded\thanks{Research was supported by the NNSF  of China (12101339, 12001526)
and Natural Science Foundation of Jiangsu Province, China (BK20210442, BK20200626).}}

\author{Chenchen Cao \\
{\small   School of Mathematics and Statistics, Ningbo University}\\
{\small Ningbo, 315211, P.R. China}\\
{\small E-mail:caochenchen@nbu.edu.cn}\\
{Zhenfeng Wu}\\
{\small School of Science,  Jiangnan University}\\
{\small Wuxi, 214122, P. R. China}\\
{\small E-mail: zhfwu@jiangnan.edu.cn}\\
{Chi Zhang\thanks{Corresponding author}}\\
{\small Department of Mathematics, China University of Mining and Technology}\\
{\small Xuzhou, 221116, P. R. China}\\
{\small E-mail: zclqq32@cumt.edu.cn}\\}

\date{}
\maketitle

\begin{abstract}
A normal subgroup $E$ of a group $G$ is said to be hypercyclically embedded in $G$ if either
$E=1$ or $E\neq 1$ and every chief factor of $G$ below $E$ is cyclic.
In this article, we present some new characterizations of a normal subgroup to be hypercyclically embedded.
Some recent results in this literature are generalized and unified.
\end{abstract}

\let\thefootnoteorig\thefootnote
\renewcommand{\thefootnote}{\empty}


 \let\thefootnote\thefootnoteorig

\section{Introduction}
In this paper, all groups considered are finite and we always use $G$ to denote a finite group.
$\mathbb{P}$ is the set of all primes and $n$ is an integer.
The symbol $\pi(n)$ denotes the set of all primes dividing $n$; in particular, $\pi(G)=\pi(|G|)$ is
the set of all primes dividing the order $|G|$ of $G$.

In order to convenience, we first introduce some notations and terminologies which often used in this paper.
Let $\sigma=\{\sigma_{i}|i\in I\}$ be some partition of $\mathbb{P}$ and $\s(G)=\{\s_i|\s_i\cap \pi(G)\neq \emptyset\}$.
$\Pi$ is a non-empty subset of the set $\s$ and $\Pi'$ denotes $\s \backslash \Pi$.
Following \cite{AN1,WS1,zs}, $G$ is said to be \emph{$\sigma$-primary} if $|\sigma(G)| \leq 1$.
An integer $n$ is said to be a \emph{$\Pi$-number} if $\pi(n)\subseteq \bigcup_{\s_i\in \Pi}\s_i$;
 a subgroup $H$ of $G$ is said to be a \emph{$\Pi$-subgroup} of $G$
if $|H|$ is a $\Pi$-number.
A set $\mathcal{H}$ of subgroups of $G$ is said to be a
\emph{complete
Hall $\s$-set} of $G$ if every non-identity element in $\mathcal{H}$
is a Hall $\s_i$-subgroup of $G$ for some $i$ and $\mathcal{H}$ contains
exactly one Hall $\s_i$-subgroup for every $\s_i\in \s (G)$.
$G$ is said to be \emph{$\s$-full} if $G$ possesses a complete Hall $\s$-set; a
\emph{$\s$-full group of Sylow type} if every subgroup of $G$ is a $D_{\s_i}$-group for
all $\s_i\in \s(G)$.
A subgroup $A$ of $G$ is said to be \emph{$\s$-subnormal} in $G$ if there exists a subgroup
chain $A=A_0\leq A_1\leq \cdots \leq A_n=G$ such that either $A_{i-1}$ is normal
in $A_i$ or $A_i/(A_{i-1})_{A_i}$ is $\s$-primary for all $i\in \{1,\cdots,n\}$.

Moreover, a subgroup $H$ of $G$ is said to be
\emph{$\sigma$-permutable} in $G$ \cite{AN1} if
$G$ possesses a complete Hall $\sigma$-set $\mathcal {H}$
such that $HA^{x}=A^{x}H$ for all $A\in \mathcal {H}$ and all $x\in G$;
\emph{$m$-$\sigma$-permutable} \cite{wei} if
$H=\langle A, B\rangle$ for some modular subgroup (in the sense of Kurosh \cite[p.43]{Sch}) $A$ and $\sigma$-permutable subgroup $B$ of $G$;
 \emph{weakly $m$-$\sigma$-permutable} in $G$ \cite{wei} if
there exists a $\sigma$-subnormal subgroup $T$ of $G$
such that $G=HT$ and $H\cap T\leq S\leq H$ for some $m$-$\s$-permutable subgroup $S$ of $G$;
\emph{$\mathcal {H}$-permutable}\cite{H-permute}
if $HA=AH$ for some complete Hall $\s$-set $\mathcal{H}$ of $G$ and all $A\in \mathcal{H}$.

Recall that a normal subgroup $E$ of $G$ is said to be \emph{hypercyclically embedded} in $G$ (see \cite[p.217]{Sch})
and denoted by $E\leq Z_{\mathfrak{U}}(G)$
if either $E=1$ or $E\neq 1$ and every chief factor of $G$ below $E$ is cyclic.
It is known that hypercyclically embedded subgroups play an
important role in the theory of finite groups $($see the books \cite{AB2, W, Sch, MW}$)$
and the permutability and modularity of subgroups have an essential influence on
a normal subgroup to be hypercyclically embedded(see, for example, \cite{MSC,Sch}).
Following this idea, the researchers studied the generalized permutability and modularity of subgroups
and obtained some interesting results under which a normal subgroup is hypercyclically embedded
(see the recent papers \cite{AN1, zs,H-permute,AN5,BL,WS2,wei,MTH,ZWG}).
In particular,
M. T. Hussain et. al \cite{MTH} prove the following result recently.

\begin{Theorem}\label{minc1}
Let $G$ be a $\s$-full group of Sylow type, $E$ a normal subgroup of $G$ and
$\mathcal {H}$ be a complete Hall $\sigma$-set of $G$
such that every member of $\mathcal{H}$ is supersoluble.
If every cyclic subgroup $H$ of any non-cyclic Hall $\s_i$-subgroup of $E$ of prime order and order $4$ $($if the Sylow $2$-subgroup of $E$ is
non-abelian and $H\nleq Z_{\infty}(G))$ is weakly $m$-$\s$-permutable in $G$ for all $\s_i\in \s(E)$,
then $E\leq Z_{\mathfrak{U}}(G)$.
\end{Theorem}

As a generalization of $\sigma$-permutable subgroups, Guo et al. \cite{H-permute}
proposed the concept of $\mathcal {H}$-permutable subgroups.
Along this way, can we generalize Theorem \ref{minc1} to more general case
by combining $\mathcal {H}$-permutable subgroups and weakly $m$-$\s$-permutable subgroups?

In this paper,
we give the definition of weakly $m$-$\mathcal {H}$-permutable subgroups
(which cover both $\mathcal {H}$-permutable subgroups and weakly $m$-$\s$-permutable subgroups, see examples below)
and obtain some new characterizations of a normal subgroup to be hypercyclically embedded
by weakly $m$-$\mathcal {H}$-permutability of subgroups.

\begin{Definition}
\rm Let $\mathcal{H}$ be a complete Hall $\s$-set of $G$. We say that a subgroup $H$ of $G$ is:

$(i)$ \emph{$m$-$\mathcal {H}$-permutable}
if $H=\langle A,B\rangle$ for some modular subgroup $A$ and $\mathcal{H}$-permutable subgroup $B$ of $G$;

$(ii)$ \emph{weakly $m$-$\mathcal {H}$-permutable} in $G$
if there exists a $\sigma$-subnormal subgroup $T$ of $G$
such that $G=HT$ and $H\cap T\leq S \leq H$ for some $m$-$\mathcal {H}$-permutable subgroup $S$ of $G$.
\end{Definition}

Before continuing, let's consider the following fact:
it is clear that every $\mathcal {H}$-permutable subgroup, every $m$-$\mathcal {H}$-permutable and every modular subgroup of $G$
are weakly $m$-$\mathcal {H}$-permutable in $G$.
But the following examples show that the converse is not true in general.

\begin{Example}\label{e1}
\rm
(i) Let $G=A_4$, where $A_4$ is
the alternating group of degree $4$.
Let $\sigma=\{\sigma_{1},\sigma_{2}\}$,
where $\sigma_{1}=\{2,3\}$ and $\sigma_{2}=\{2,3\}^{'}$.
Then, clearly, $\mathcal {H}=\{1,G\}$ is a complete Hall $\sigma$-set of $G$ and every subgroup of $G$ is
$\mathcal {H}$-permutable.
It follows that every subgroup of $G$ is weakly $m$-$\mathcal {H}$-permutable in $G$.
Let $H$ be a subgroup of order $2$ of $G$.
In view of \cite[Theorem 5.1.14]{Sch}, $H$ is not modular in $G$.

(ii) Let $p, q,r, t$ be distinct primes, where $q$ divides $p-1$ and $t$ divides $r -1$.
Let $\sigma=\{\sigma_{1},\sigma_{2},\sigma_{3}\} $,
where $\sigma_{1}=\{p,q\}$, $\sigma_{2}=\{r\}$ and $\sigma_{3}=\{p,q,r\}^{'}$.
Let $V = Q \rtimes C_p$,
where $Q$ is a simple $F_qC_p$-module which is faithful for $C_p$. Let $C_r  \rtimes C_t$ be a non-abelian group of order $rt$, $A = C_t$
and
$G = (Q  \rtimes C_p) \times (C_r  \rtimes C_t)$. Let $T$ be a subgroup of $G$ of order $t$ such that $T\neq A$.
Then $\mathcal {H}=\{1,V,C_r,T \}$ is a complete Hall $\s$-set of $G$. Let $B$ be a subgroup of order $q$ in $Q$. Then $B < Q$ since $p > q$,
and so $B$ is $\mathcal {H}$-permutable. Let
$H = \langle A, B\rangle=A\times B$. In view of \cite[Theorem 5.1.9]{Sch}, we have that $A$ is modular in $G$. Hence $H$ is $m$-$\mathcal {H}$-permutable,
and so weakly $m$-$\mathcal {H}$-permutable in $G$.
Assume that $H$ is $\mathcal {H}$-permutable. Then $HT=TH$, which shows that $HT$ is a subgroup of $G$.
Since $A$ and $T$ are both Sylow $t$-subgroups of $HT$, there exists an element $g=th \in HT$ such that
$A=T^{th}=T^h$, where $t\in T$ and $h\in H$.
But since $h\in H=A\times B$, we have that $h=ab$ for some elements $a\in A$ and $b\in B$.
It follows from $A=T^{h}$ that $T=A^{h^{-1}}=A^{b^{-1}a^{-1}}=A$, contrary to the choice of $T$.
Hence $H$ is not $\mathcal {H}$-permutable.
\end{Example}

Note that all $\sigma$-permutable subgroups and all weakly
$m$-$\sigma$-permutable subgroups of $G$ are weakly $m$-$\mathcal {H}$-permutable in $G$ for some
complete Hall $\s$-set $\mathcal{H}$ of $G$.
But the following example shows that weakly $m$-$\mathcal {H}$-permutable subgroups of $G$ are not necessarily
weakly $m$-$\s$-permutable in $G$.

\begin{Example}
\rm
Let $G=A_5$, where $A_5$ is
the alternating group of degree $5$.
Let $H$ be a subgroup of $G$ of order $5$
and $K\cong A_4$ be a subgroup of $G$,
where $A_4$ is
the alternating group of degree $4$.
Let $\sigma=\{\sigma_{1},\sigma_{2}\}$,
where $\sigma_{1}=\{2,3\}$ and $\sigma_{2}=\{2,3\}^{'}$.
Clearly, $\mathcal{H}=\{H,K\}$ is a Hall $\s$-set of $G$
and $H$ is $\mathcal{H}$-permutable, and so
$H$ is weakly $m$-$\mathcal{H}$-permutable in $G$.
Suppose that $H$ is weakly $m$-$\s$-permutable in $G$.
Then there exists a
$\s$-subnormal subgroup $T$
such that $G=HT$ and $H\cap T\leq S\leq H$ for some $m$-$\s$-permutable subgroup $S$ of $G$.
Since $G$ is a simple group and $|\sigma(G)|>1$, we obtain that $1$ and $G$ are the only two $\s$-subnormal subgroups of $G$.
It follows that $T=G$, and so $H=S$.  Let $S=\langle A,B \rangle$, where $A$ is modular and $B$ is $\sigma$-permutable in $G$.
If $A\neq 1$, then every chief factor of $G$ between $A^G$ and $A_G$ is cyclic (see Lemma \ref{cyclic} below), a contradiction.
Hence $A= 1$, and so $H=S=B$ is $\sigma$-permutable in $G$. This implies that $HH^g=H^gH$ for any element $g\in G$, which is impossible.
Hence $H$ is not weakly $m$-$\s$-permutable in $G$.
\end{Example}

In fact, our main goal here is to prove the following more general result.

\begin{Theorem}\label{minc}
Let $G$ be a $\s$-full group of Sylow type, $E$ a normal subgroup of $G$ and
$\mathcal {H}$ be a complete Hall $\sigma$-set of $G$
such that every member of $\mathcal{H}$ is supersoluble.
If every cyclic subgroup $H$ of any non-cyclic Hall $\s_i$-subgroup of $E$ of prime order and order $4$ $($if the Sylow $2$-subgroup of $E$ is
non-abelian and $H\nleq Z_{\infty}(G))$ is weakly $m$-$\mathcal{H}$-permutable in $G$ for all $\s_i\in \s(E)$,
then $E\leq Z_{\mathfrak{U}}(G)$.
\end{Theorem}

\begin{Remark}
{\rm In the classical case when $\s=\{\{2\},\{3\},\cdots\}$, $\s$-permutable subgroups
are also called
$s$-permutable subgroups.
A subgroup $H$ of $G$ is said to be $c$-normal in $G$ \cite{WY1} if $G$ has a normal subgroup $T$
such that $G=HT$ and $H\cap T\leq H_G$, where $H_G$ denote the maximal normal subgroup of $G$ contained in $H$.
Hence by Theorems \ref{minc}, we may directly obtain many known results, for example,
\cite[Theorem 1.3 and Corollary 1.6]{MTH}, \cite[Theorem 4.2]{WY1}, \cite[Corollary 3.6]{Shaa},
\cite[Theorem 3.1]{SS} and \cite[Theorem 3]{Bu}.}
\end{Remark}

This paper is organized as follows. In Section 2, we cite some known results which are useful
in our proofs and prove some basic properties of weakly $m$-$\mathcal{H}$-permutable
subgroups. In section 3, we give the proofs of Theorem \ref{minc}.
In section $4$, we give some more applications of our
result.

\section{Preliminaries}

Following \cite{AN1,WS1}, we use $O^{\Pi}(G)$ to denote the subgroup of $G$ generated by all its $\Pi'$-subgroups.
Instead of $O^{\{\s_i\}}(G)$, we write $O^{\s_i}(G)$.

\begin{Lemma}\label{subnormal}{\rm\cite[Lemma 2.6]{AN1}}
Let $A,K$ and $N$ be subgroups of $G$. Suppose that $A$ is $\sigma$-subnormal in $G$ and $N$ is normal in $G$.

$(1)$ $A\cap K$ is $\sigma$-subnormal in $K$.

$(2)$ $AN/N$ is $\sigma$-subnormal in $G/N$.

$(3)$ If $|G:A|$ is a $\Pi$-number,
then $O^{\Pi}(A)=O^{\Pi}(G)$.
\end{Lemma}

Suppose that $G$ possesses a complete Hall $\s$-set $\mathcal{H}=\{H_1,\cdots,H_t\}$. Following \cite{H-permute}, for any
subgroup $H$ (resp. normal subgroup $N$) of $G$ we write $H\cap \mathcal{H}$ (resp. $\mathcal{H}N/N$)
to denote the set $\{H\cap H_1,\cdots,H\cap H_t\}$ (resp. $\{H_1N/N,\cdots,H_tN/N\}$).
If $H\cap \mathcal{H}$ is a complete Hall $\s$-set of $H$, then we say that $\mathcal{H}$ \emph{reduces into} $H$.

\begin{Lemma}\label{permutable}
Let $\mathcal{H}$ be a complete Hall $\s$-set of $G$ and $R$ a normal subgroup of $G$.
Suppose that $H\leq E\leq G$ and $H$ is $m$-$\mathcal{H}$-permutable in $G$.

$(1)$ $HR/R$ is $m$-$($$\mathcal{H}R/R$$)$-permutable in $G/R$ .

$(2)$ If $\mathcal{H}$ reduces into $E$, then $H$ is $m$-$($$E\cap \mathcal{H}$$)$-permutable in $E$.

$(3)$ If $E$ is normal in $G$, then $\mathcal{H}$ reduces into $E$ and so $H$ is $m$-$($$E\cap \mathcal{H}$$)$-permutable in $E$.

$(4)$ The group $\langle H,K\rangle$ is $m$-$\mathcal{H}$-permutable in $G$ for any $m$-$\mathcal{H}$-permutable subgroup $K$ of $G$.

\end{Lemma}

\demo
By the hypothesis, we may assume that $H=\langle A,B \rangle$, where $A$ is modular and $B$ is $\mathcal{H}$-permutable.

$(1)$ In view of \cite[p.201, Property (3)]{Sch}, we have that $AR/R$ is modular in $G/R$.
It follows from $B$ is $\mathcal{H}$-permutable that $BN/N$ is $(\mathcal{H}R/R)$-permutable
by \cite[Lemma 2.1(1)]{H-permute}.
Moreover, it is clear that $HR/R=\langle AR/R,BR/R \rangle$. Hence $HR/R$ is $m$-$($$\mathcal{H}R/R$$)$-permutable in $G/R$.

$(2)$ \cite[p.201, Property (2)]{Sch} shows that $A$ is modular in $E$.
Besides, since $\mathcal{H}$ reduces into $E$, we obtain that $B$ is $($$E\cap \mathcal{H}$$)$-permutable by \cite[Lemma 2.1(2)]{H-permute}.
This shows that $H$ is $m$-$($$E\cap \mathcal{H}$$)$-permutable in $E$.

$(3)$ Lemma 3.2(c) in \cite[Chapter I]{Doerk} implies that $\mathcal{H}$ reduces into $E$, and so $H$ is $m$-$($$E\cap \mathcal{H}$$)$-permutable in $E$ by $(2)$.

$(4)$ Assume that $K=\langle C,D \rangle$, where $C$ is modular and $D$ is $\mathcal{H}$-permutable.
Then, clearly, $\langle H, K\rangle=\langle \langle A,C\rangle,\langle B,D\rangle\rangle$. Then
$\langle A,C\rangle$ is modular in $G$ by \cite[p.201, Property (5)]{Sch}.
Besides, since $\langle B,D\rangle$ is $\mathcal{H}$-permutable
$($see \rm\cite[Chapter A, 1.6(a)]{Doerk}$)$,
we obtain that $\langle H,K\rangle$ is $m$-$\mathcal{H}$-permutable in $G$.
\qed

\begin{Lemma}\label{main}
Let $\mathcal{H}$ be a complete Hall $\s$-set of $G$.
Suppose that $H\leq E\leq G$ and $H$ is weakly $m$-$\mathcal{H}$-permutable in $G$.

$(1)$ If $\mathcal{H}$ reduces into $E$,
then $H$ is weakly $m$-$($$E\cap \mathcal{H}$$)$-permutable in $E$.

$(2)$ If $E$ is normal in $G$, then $\mathcal{H}$ reduces into $E$
and so $H$ is weakly $m$-$($$E\cap \mathcal{H}$$)$-permutable in $E$.

$(3)$ If $N$ is a normal subgroup of $G$ with $N\leq H$,
then $H/N$ is weakly $m$-$($$\mathcal{H}N/N$$)$-permutable in $G/N$.

$(4)$ If $N$ is a normal subgroup of $G$ with $(|H|,|N|)=1$,
then $HN/N$ is weakly $m$-$($$\mathcal{H}N/N$$)$-permutable in $G/N$.
\end{Lemma}

\demo
By the hypothesis, there exists a $\sigma$-subnormal subgroup $T$ of $G$
such that $G=HT$ and $H\cap T\leq S\leq H$ for some $m$-$\mathcal{H}$-permutable subgroup $S$ of $G$.

$(1)$ Since $H\leq E$,
we have that $E=E\cap HT=H(E\cap T)$.
By Lemma \ref{subnormal}(1), $E\cap T$ is $\sigma$-subnormal in $E$.
Moreover, since $\mathcal{H}$ reduces into $E$, $S$ is $m$-$\mathcal{H}$-permutable in $E$ by Lemma \ref{permutable}(2).
Note that
$H\cap(E\cap T)=H\cap T\leq S\leq H$.
This shows that $H$ is weakly $m$-$(E\cap\mathcal{\mathcal{H}})$-permutable in $E$.

$(2)$ Lemma $2.3(3)$ shows that $\mathcal{H}$ reduces into $E$, and so
$H$ is weakly $m$-$(E\cap\mathcal{\mathcal{H}})$-permutable in $E$ by $(1)$.

$(3)$ Lemma $3.2(b)$ in \cite[Chapter I]{Doerk} implies that $\mathcal{H}N/N$ is a complete Hall $\s$-set
of $G/N$.
Clearly, $G/N=(H/N)(TN/N)$.
By Lemma \ref{permutable}(1), $SN/N$ is $m$-$($$\mathcal{H}N/N$$)$-permutable in $G/N$.
Besides, since $H/N\cap TN/N=(H\cap T)N/N\leq SN/N\leq H/N$ and $TN/N$ is $\sigma$-subnormal in $G/N$ by
Lemma \ref{subnormal}(2), we get that $H/N$ is weakly $m$-$(\mathcal{H}N/N)$-permutable in $G/N$.

$(4)$
It is clear that $\mathcal{H}N/N$ is a complete Hall $\s$-set
of $G/N$.
Since $(|H|,|N|)=1$,
we have that $(|G:T|,|N|)=1$.
It follows from $|G:T|=|G:NT||NT:T|=|G:NT||N:N\cap T|$ that $N\cap T=N$,
and so $N\leq T$. Hence $G/N=(HN/N)(T/N)$.
By Lemma \ref{subnormal}(2), we know that $T/N$ is $\sigma$-subnormal in $G/N$.
It follows from Lemma \ref{permutable}(1) that $SN/N$ is $m$-$($$\mathcal{H}N/N$$)$-permutable in $G/N$.
Moreover, since $(HN/N)\cap (T/N)=(H\cap T)N/N\leq SN/N\leq HN/N$,
we obtain that $HN/N$ is weakly $m$-$(\mathcal{H}N/N)$-permutable in $G/N$.
\qed

\begin{Lemma}{\rm\cite[Lemma 5]{BNVS}}\label{cap}
Let $H,K$ and $N$ be pairwise permutable subgroups of $G$, and suppose that $H$ is a Hall subgroup of $G$. Then $N\cap HK=(N\cap H)(N\cap K)$.
\end{Lemma}

\begin{Lemma}\label{cyclic}{\rm\cite[Theorem 5.2.5]{Sch}}
If $H$ is modular in $G$, then every chief factor of $G$ between $H^G$ and $H_G$ is cyclic.
\end{Lemma}

Let $P$ be a $p$-group. If $P$ is not a non-abelian $2$-group, then we use $\Omega(P)$ to denote $\Omega_{1}(P)$. Otherwise, $\Omega(P)=\Omega_2(P)$.

\begin{Lemma}{\rm\cite[Lemma 4.3]{GSO}}\label{Tc}
Let $C$ be a Thompson critical subgroup of a $p$-group $P$ $($see {\rm \cite[p.185]{Gfini}}$)$.

$(1)$ If $p$ is odd, then the exponent of $\Omega(C)$ is $p$.

$(2)$ If $P$ is a non-abelian $2$-group, then the exponent of $\Omega(C)$ is $4$.
\end{Lemma}

\begin{Lemma}{\rm\cite[Lemma 4.4]{GSO}}\label{super}
Let $P$ be a normal $p$-subgroup of $G$ and $C$ be a Thompson critical subgroup of $P$.
If $\Omega(C)$ is hypercyclically embedded in $G$, then $P$ is hypercyclically embedded in $G$.
\end{Lemma}

\section{Proof of Theorem \ref{minc}}

In order to prove our main result, we may first prove the following Propositions.

\begin{Proposition}\label{prop1}
Let $G$ be a $\sigma$-full group of Sylow type
and $\mathcal {H}$ a complete Hall $\sigma$-set of $G$
such that every member of $\mathcal{H}$ is supersoluble.
Suppose that the smallest prime dividing $|G|$ is $p$.
If every cyclic subgroup $H$ of order $p$ and order $4$ $($if $p=2$, a Sylow $2$-subgroup of $G$ is
non-abelian and $H\nleq Z_{\infty}(G))$ is weakly $m$-$\mathcal{H}$-permutable in $G$,
then $G$ is $p$-nilpotent.
\end{Proposition}

\demo
Assume that this is false and let $G$ be a counterexample. Then $G$ is not $p$-nilpotent. Hence by {\rm \cite[IV, Theorem ~5.4]{HU}}
and {\rm \cite[Theorem ~3.4.11]{W1}}, $G$ has a $p$-closed Schmidt subgroup
$I=P\rtimes Q$, where $P$ is a Sylow $p$-subgroup of $I$ of exponent $p$ or $4$
$($if $P$ is a non-abelian $2$-group$)$, $Q$ is a Sylow $q$-subgroup for some prime $q\neq p$,
$P/\Phi(P)$ is an $I$-chief factor, $Z_{\infty}(I)=\Phi(I)$ and $\Phi(I)\cap P=\Phi(P)$.

Without loss of generality, we may assume that $p\in \sigma_1$ and let $H_1\in \mathcal{H}$
be a Hall $\sigma_1$-subgroup of $G$.
We claim that $|P:\Phi(P)|=p$.
If $q\in \sigma_1$, then $I$ is a $\s_1$-subgroup of $G$, and so $I\leq H_1^g$
for some $g\in G$ by $G$ is a $\s$-full group of Sylow type.
Since $H_1$ is supersoluble and $P/\Phi(P)$ is an $I$-chief factor, we have that $|P:\Phi(P)|=p$.
Now, we consider that $q\notin \sigma_1$.
Assume that $X/\Phi(P)$ is a minimal subgroup of $P/\Phi(P)$.
Let $x\in X\backslash \Phi(P)$ and $L=\langle x\rangle$.
Then $X=L\Phi(P)$ and $|L|=p$ or $4$ $($if $P$ is a non-abelian $2$-group$)$.
If $L\leq Z_{\infty}(G)$, then $L\leq Z_{\infty}(I)\cap P=\Phi(I)\cap P=\Phi(P)$, which contradicts the choice of $L$.
Hence $L\nleq Z_{\infty}(G)$.
Then by the hypothesis, $L$ is weakly $m$-$\mathcal{H}$-permutable in $G$.
That is, there exists a $\s$-subnormal subgroup $T$ of $G$ such that $G=LT$ and $L\cap T\leq S\leq L$ for some $m$-$\mathcal {H}$-permutable subgroup $S$ of $G$.
It follows that $I=I\cap LT=L(I\cap T)$ and $(I\cap T)$ is $\s$-subnormal in $I$ by Lemma \ref{subnormal}(1).
We write $T_1=I\cap T$.
Since $T_1$ is $\s$-subnormal in $I$ and $|I:T_1|$ is a $\s_1$-number,
we obtain that $O^{\s_1}(T_1)=O^{\s_1}(I)$ by Lemma \ref{subnormal}$(3)$.
But $q\notin \sigma_1$ and $Q\leq O^{\s_1}(I)$, so $Q^I\leq O^{\s_1}(I)$.
If $Q^I\neq I$, then $Q$ $char$ $Q^I \trianglelefteq I$ since $I$ is a Schmidt group, and so $Q \trianglelefteq I$, a contradiction.
Hence $I=Q^I\leq O^{\s_1}(I)$. It follows from $O^{\s_1}(T_1)=O^{\s_1}(I)$ that $T_1=I$.
Therefore $L=L\cap T_1\leq L\cap T\leq S\leq L$, and so $L=S$.
Let $S=\langle A,B \rangle$, where $A$ is modular and $B$ is $\mathcal{H}$-permutable.
In view of \cite[p.201, Property (2)]{Sch}, we have that $A$ is modular in $I$, and so the chief factor of $I$ between $A^I$ and $A_I$ is cyclic by Lemma \ref{cyclic}.
Hence the chief factor of $I/\Phi(P)$ between $A^I\Phi(P)/ \Phi(P)$ and $A_I\Phi(P)/\Phi(P)$ is cyclic.
Moreover, since $A\Phi(P)/ \Phi(P)\leq L\Phi(P)/\Phi(P)$ and $L\Phi(P)/\Phi(P)=X/\Phi(P)$ is a minimal subgroup of $P/\Phi(P)$,
we obtain that $A\Phi(P)/ \Phi(P)=L\Phi(P)/ \Phi(P)$ or $A\Phi(P)/ \Phi(P)=\Phi(P)/ \Phi(P)$.
In the former case, we get that $A^I\Phi(P)/ \Phi(P)= (A\Phi(P))^I/ \Phi(P) =   (L\Phi(P))^I/ \Phi(P)=P/\Phi(P)$ and
$A_I\Phi(P)/ \Phi(P)= (A\Phi(P))_I/ \Phi(P) =   (L\Phi(P))_I/ \Phi(P)=\Phi(P)/\Phi(P)$
since $P/\Phi(P)$ is an $I$-chief factor. It follows that $P/\Phi(P)$ is cyclic, and so $|P:\Phi(P)|=p$, as desired.
In the later case, we obtain that $A\leq \Phi(P)$, and so $A< S=L$ since $L\nleq \Phi(P)$.
It follows from $L$ is cyclic group of order $p$ or $4$ that $A\leq \Phi(L)$.
Hence $L=\langle A,B\rangle =B$,
which shows that $L$ is $\mathcal{H}$-permutable.
 That is,
$LD=DL$ for all $D\in \mathcal{H}$. Assume that $D$ is a
Hall $\s_j$-subgroup of $G$, where $q\in \s_j$. Then
$Q^g\leq D$ for some $g\in G$.
Next, we consider the group $I^g=P^g\rtimes Q^g$ and the subgroup $L^g$ of $I^g$.
By a similar argument as above, we have that $|P/\Phi(P)|=|P^g/\Phi(P)^g|=p$ or $L^gD=D L^g$.
Without loss of generality, we only need to consider the case when $L^g  D =D L^g$.
Since $L$ is cyclic, we have that $L^g  D =D L^g$ is $p$-nilpotent
by {\rm \cite[IV, Theorem 2.8]{HU}}. Hence $L^g\leq N_G(D)$.
We choose any element $s\in Q$.
Then $[s,x]=s^{-1}x^{-1}sx=s^{-1}s^{x}=s^{-1}s^{gg^{-1}xgg^{-1}}=s^{-1}s^{gx^gg^{-1}}\in  D^{g^{-1}}$.
On the other hand,
we have that $[s,x]=s^{-1}x^{-1}sx=(x^{-1})^{s}x\in P$.
Notice that $(|P|, | D^{g^{-1}}|)=1$, we obtain that $[s,x]=1$, and so $L\leq C_I(Q)$.
It follows from $L\Phi(P)/ \Phi(P)\unlhd P/\Phi(P)$ that $L\Phi(P)/ \Phi(P)\unlhd I/\Phi(P)$.
Hence $L\Phi(P)/ \Phi(P)=P/\Phi(P)$ for $P/\Phi(P)$ is an $I$-chief factor, which shows that $|P/\Phi(P)|=p$.

Therefore $P$ is cyclic and so $P$ is of exponent $p$ by above.
This implies that $P$ is a group of order $p$.
Since $N_I(P)/C_I(P)\lesssim Aut(P)$ is a group of order $p-1$ and $p$ is the smallest prime dividing $|E|$,
it follows that $C_I(P)=N_I(P)=I$.
Thus $Q\trianglelefteq I$. This contradiction completes the proof.
\qed

\begin{Proposition}\label{mincyc}
Let $G$ be a $\s$-full group of Sylow type and $\mH=\{H_1,\cdots,H_t\}$ a complete Hall $\s$-set of $G$ such that $H_i$ is a
supersoluble $\s_i$-subgroup for all $i\in \{1,\cdots,t\}$.
Let $P$ be a normal $p$-group of $G$ and $P\leq H_j$ for some $j$.
If every cyclic subgroup $H$ of $P$ of prime order and order $4$ $($if $P$ is a
non-abelian $2$-group and $H\nleq Z_{\infty}(G))$ is weakly $m$-$\mathcal{H}$-permutable in $G$,
then $P \leq Z_{\mathfrak{U}}(G)$.
\end{Proposition}

\demo
Assume that this is false and let $(G,P)$ be a counterexample with $|G|+|P|$ minimal.
Without loss of generality, we may assume that $j=1$.

$(1)$ {\sl Let $P/N$ be a chief factor of $G$.
Then $N \leq Z_{\mathfrak{U}}(G)$.
Hence $N$ is the unique normal subgroup of $G$ satisfying that
$P/N$ is a chief factor of $G$, and $|P/N|>p$}.

Clearly, $(G,N)$ satisfies the hypothesis of this reslut.
Hence $N \leq Z_{\mathfrak{U}}(G)$ by the choice of $(G,P)$.
Assume that $G$ has another normal subgroup $R\neq N$ such that $P/R$ is a chief factor of $G$.
Then $R \leq Z_{\mathfrak{U}}(G)$.
It follows that $P/N=RN/N \leq Z_{\mathfrak{U}}(G/N)$,
and so $P \leq Z_{\mathfrak{U}}(G)$, a contradiction.
It is clear also that $|P/N|>p$ by above analysis.

$(2)$ {\sl The exponent of $P$ is $p$ or $4$ $($if $P$ is a non-abelian $2$-group$)$}.

Let $C$ be a Thompson critical subgroup of $P$.
If $\Omega(C)< P$, then $\Omega(C)\leq N \leq Z_{\mathfrak{U}}(G)$ by Claim $(1)$.
It follows that $P \leq Z_{\mathfrak{U}}(G)$ by Lemma $\ref{super}$, a contradiction.
Hence we have $\Omega(C)= P$.
By Lemma $\ref{Tc}$, we have that the exponent of $P$ is $p$ or $4$ $($if $P$ is a non-abelian $2$-group$)$.

$(3)$ {\sl Final contradiction}.

Since $H_1/N$ is supersoluble and $|P/N|>p$, there exists a minimal normal subgroup $L/N$ of $H_1/N$ such that $1\neq L/N<P/N$ and $L/N$ is cyclic.
Let $x\in L\setminus N$ and $H=\langle x\rangle$. Then $L=HN$ and $|H|=p$ or $4$ $($if $P$ is a non-abelian $2$-group$)$ by Claim $(2)$.
Assume that $H\leq Z_{\infty}(G)$, then $L/N=HN/N\leq Z_{\infty}(G)N/N\leq Z_{\infty}(G/N)$ by {\rm \cite[Chapter 1, Theorem 2.6(d)]{W}},
which means that $Z_{\infty}(G/N)\cap P/N\neq 1$. Hence $P/N\leq Z_{\infty}(G/N)$ since $P/N$ is a chief factor of $G$.
It follows from Claim $(1)$ that $P \leq Z_{\mathfrak{U}}(G)$.
This contradiction shows that $H\nleq Z_{\infty}(G)$.
Then by the hypothesis, there exists a $\s$-subnormal subgroup $T$ of $G$ such that
$G=HT$ and $H\cap T\leq S\leq H$ for some $m$-$\mathcal {H}$-permutable subgroup $S$ of $G$.
Hence by Lemma $\ref{subnormal}(3)$, we have that $O^{\s_1}(T)=O^{\s_1}(G)$ since $|G:T|$ is a $\s_1$-number.

We claim that $O^{\s_1}(G/N) \leq N_{G/N}(L/N)$. First assume that $P\cap O^{\s_1}(G)=P$, then $H\leq P\leq O^{\s_1}(G)=O^{\s_1}(T)\leq T$.
Hence $G=HT=T$, and so $H=H\cap T\leq S\leq H$, which shows that $H=S$.
Let $S=\langle A,B \rangle$, where $A$ is modular and $B$ is $\mathcal{H}$-permutable.
Since $A\leq H<P$, we have that $A^G\leq P$ and $A_G\leq N$ by Claim (1).
If $A^G=P$, then $P/N$ is cyclic by Lemma \ref{cyclic}, contrary to Claim $(1)$.
Consequently, $A^G\leq N$ by Claim $(1)$ again.
Hence $L/N=HN/N=SN/N=\langle A,B\rangle N/N=BN/N$ is $(\mathcal{H}N/N)$-permutable by Lemma \cite[Lemma 2.1(1)]{H-permute}.
 That is,
$(L/N) (DN/N)=(DN/N)(L/N)$ for all $D\in \mathcal{H}$. Assume that $D$ is not a
Hall $\s_1$-subgroup of $G$. Then $DN/N$ is not a
Hall $\s_1$-subgroup of $G/N$
and $P/N\cap (L/N) (DN/N)=(P/N \cap L/N)(P/N\cap DN/N)=L/N$
by Lemma \ref{cap}. Hence $DN/N\leq N_G(L/N)$,
which implies that $O^{\s_1}(G/N) \leq N_{G/N}(L/N)$.
Now, assume that $P\cap O^{\s_1}(G)<P$.
Then $P\cap O^{\s_1}(G)\leq N$ by Claim $(1)$.
Hence $O^{\s_1}(G/N)\cap P/N\leq O^{\s_1}(G)N/N\cap P/N=(O^{\s_1}(G)\cap P)N/N=1$ by Dedekind Modular Law.
It follows that $O^{\s_1}(G/N) \leq C_{G/N}(P/N)$ and so $O^{\s_1}(G/N) \leq N_{G/N}(L/N)$.
Moreover, since $L/N\trianglelefteq H_1/N$, we obtain that $L/N\trianglelefteq G/N$, a contradiction.
This completes the proof of the result.
\qed

\medskip
\textbf{Proof of Theorem \ref{minc}}~~
Assume that this is false and let $(G,E)$ be a counterexample with minimal $|G|+|E|$.
Let $p$ is the smallest prime dividing $|E|$ and $P$ a Sylow $p$-subgroup of $E$.
Without loss of generality, we may assume that $P\leq H_1\cap E$, where $H_1\in \mathcal{H}$.

$(1)$ {\sl $E$ is $p$-nilpotent}.

In view of Lemma \ref{permutable}(3), we have that the group $E$ satisfies the hypothesis of Proposition \ref{prop1}.
Hence $E$ is $p$-nilpotent by Proposition \ref{prop1}.

\smallskip
By Claim (1), we may let $E_{p'}$ be the normal Hall $p'$-subgroup of $E$. Then $E_{p'}\trianglelefteq G$.
\smallskip

$(2)$ {\sl $E_{p'}=1$, and so $E$ is a $p$-group.}

Assume that $E_{p'}\neq 1$. Then we consider the quotient group $G/E_{p'}$.
It is clear that $\mathcal{H}E_{p'}/E_{p'}=\{H_1E_{p'}/E_{p'},\cdots,H_tE_{p'}/E_{p'}\}$
is a complete Hall $\s$-set of $G/E_{p'}$ and $H_iE_{p'}/E_{p'}\cong H_i/H_i\cap E_{p'}$ is supersoluble.

We claim that the hypothesis holds for $(G/E_{p'},E/E_{p'})$.
In fact, $H_iE_{p'}/E_{p'}\cap E/E_{p'}=(H_i\cap E)E_{p'}/E_{p'}=1$ for $i\in \{2,\cdots, t\}$
and $H_1E_{p'}/E_{p'}\cap E/E_{p'}=(H_1\cap E)E_{p'}/E_{p'}=E/E_{p'}$. It is trivial when $E/E_{p'}$ is cyclic.
We may, therefore, assume that $E/E_{p'}$ is non-cyclic.
Let $H/E_{p'}$ be a cyclic subgroup of $E/E_{p'}$ of order $p$ or $4$ $($if the Sylow $2$-subgroup of $E/E_{p'}$ is non-abelian
and $H/E_{p'}\nleq Z_{\infty}(G/E_{p'})$$)$.
Then $H=E_{p'}\rtimes L$,
where $L=\langle x\rangle$ is order $p$ or $4$ $($if the Sylow $2$-subgroup of $E$ is non-abelian$)$,
by Schur-Zassenhaus theorem.
If $L\leq Z_{\infty}(G)$, then $H/E_{p'}=LE_{p'}/E_{p'}\leq Z_{\infty}(G)E_{p'}/E_{p'}\leq Z_{\infty}(G/E_{p'})$
by {\rm \cite[Chapter 1, Theorem 2.6(d)]{W}}.
Now assume that $L\nleq Z_{\infty}(G)$.
Then by Lemma \ref{main}(4), we know that the hypothesis holds for $(G/E_{p'},E/E_{p'})$.
Thus $E/E_{p'} \leq Z_{\mathfrak{U}}(G/E_{p'})$.
It is clear that the hypothesis holds for $(G,E_{p'})$.
Hence $E_{p'} \leq Z_{\mathfrak{U}}(G)$ by the choice of $(G,E)$.
Therefore $E \leq Z_{\mathfrak{U}}(G)$, a contradiction. Hence we have $(2)$.

$(3)$ {\sl Final contradiction}.

It is clear that $E$ is not cyclic, then $(G,E)$ satisfies the hypothesis of Proposition $\ref{mincyc}$ by the hypothesis and Claim (2).
Hence $E \leq Z_{\mathfrak{U}}(G)$ by Proposition $\ref{mincyc}$.
This finishes the proof.
\qed

\section{Some applications of the results}

Recall a class of groups is a \emph{formation}
if for every group $G$, every homomorphic image of $G/G^{\frak {F}}$
 belongs to $\frak {F}$.
A formation
$\mathfrak{F}$ is said to be \emph{saturated} if $G/\Phi(G)\in \mathfrak{F}$ implies
that $G\in \mathfrak{F}$.

In view of \cite[Theorem C]{AN4} and \cite[Lemma 2.16]{AN3}, we can directly get the following corollaries from Theorem \ref{minc}.

\begin{Corollary} \label{cor3}
Let $\mathfrak{F}$ be a saturated formation containing all supersoluble groups
and $E$ a normal subgroup of $G$ with $G/E\in \mathfrak{F}$.
Suppose that $G$ is a $\sigma$-full group of Sylow type
and $\mathcal {H}$ is a complete Hall $\sigma$-set of $G$
such that every member of $\mathcal{H}$ is supersoluble.
If every cyclic subgroup $H$ of any non-cyclic Hall $\s_i$-subgroup of $F^{*}(E)$ of prime order and order $4$ $($if the Sylow $2$-subgroup of $E$ is
non-abelian and $H\nleq Z_{\infty}(G))$ is weakly $m$-$\mathcal{H}$-permutable in $G$ for all $\s_i\in \s(F^{*}(E))$,
then $G\in \mathfrak{F}$.
\end{Corollary}

In the corollaries, $F^{*}(E)$ is the generalized Fitting subgroup of $E$, that is, the largest
normal quasinilpotent subgroup of $E$ (see \cite[Chapter X]{HuNB}).

Theorem \ref{minc} and Corollary \ref{cor3} also cover many other results, for example,

\begin{Corollary}{\rm\cite[Theorem 3.4]{BaW}}
Let $\mathfrak{F}$ be a saturated formation containing all supersoluble groups and $G$ a group.
If every cyclic subgroup of $G^{\mathfrak{F}}$ of prime order and order $4$ is $c$-normal in $G$,
then $G\in \mathfrak{F}$.
\end{Corollary}

\begin{Corollary}{\rm\cite[Corollary 3.5]{BaW}}
Let $G$ be a group.
If every cyclic subgroup of $G$ of prime order and order $4$ is $c$-normal in $G$,
then $G$ is supersoluble.
\end{Corollary}

\begin{Corollary}{\rm\cite[Theorem 2]{BP}}
Let $\mathfrak{F}$ be a saturated formation containing all supersoluble groups and $E$ a normal subgroup of $G$ with $G/E\in \mathfrak{F}$.
If every cyclic subgroup of $E$ of prime order and order $4$ is permutable in $G$,
then $G\in \mathfrak{F}$.
\end{Corollary}

\begin{Corollary}{\rm\cite[Theorem 5]{BP}}
Let $\mathfrak{F}$ be a saturated formation containing all supersoluble groups and $E$ a normal subgroup of $G$ with $G/E\in \mathfrak{F}$.
Suppose that $G$ has an abelian Sylow $2$-subgroup.
If every cyclic subgroup of $E$ of prime order is permutable in $G$,
then $G\in \mathfrak{F}$.
\end{Corollary}

\begin{Corollary}{\rm\cite[Theorem]{ACs}}
Let $\mathfrak{F}$ be a saturated formation containing all supersoluble groups and $G$ a group.
Suppose that $G$ has a soluble normal subgroup $E$ with $G/E\in \mathfrak{F}$.
If every cyclic subgroup of $F(E)$ of prime order and order $4$ is $s$-permutable in $G$,
then $G\in \mathfrak{F}$.
\end{Corollary}

\begin{Corollary}{\rm\cite[Corollary 1]{ACs}}
Suppose that $G$ has a soluble normal subgroup $E$ with $G/E$ is supersoluble.
If every cyclic subgroup of $F(E)$ of prime order and order $4$ is $s$-permutable in $G$,
then $G$ is supersoluble.
\end{Corollary}

\begin{Corollary}{\rm\cite[Corollary 2]{ACs}}
Let $G$ be a soluble group.
If every cyclic subgroup of $F(E)$ of prime order and order $4$ is $s$-permutable in $G$,
then $G$ is supersoluble.
\end{Corollary}

\begin{Corollary}{\rm\cite[Theorem 3.1]{LiWa}}
Let $N$ be a normal subgroup with $G/N$ is supersoluble.
If every cyclic subgroup of $F^{*}(N)$ of prime order and order $4$ is $s$-permutable in $G$,
then $G$ is supersoluble.
\end{Corollary}

\begin{Corollary}{\rm\cite[Theorem 3.3]{LiWa}}
Let $\mathfrak{F}$ be a saturated formation containing all supersoluble groups and $G$ a group.
Suppose that $G$ has a soluble normal subgroup $E$ with $G/E\in \mathfrak{F}$.
If every cyclic subgroup of $F^{*}(E)$ of prime order and order $4$ is $s$-permutable in $G$,
then $G\in \mathfrak{F}$.
\end{Corollary}

\begin{Corollary}{\rm\cite[Theorem 2]{WH}}
Let $\mathfrak{F}$ be a saturated formation containing all supersoluble groups and $G$ a group.
Suppose that $G$ has a soluble normal subgroup $E$ with $G/E\in \mathfrak{F}$.
If every cyclic subgroup of $F(E)$ of prime order and order $4$ is $c$-normal in $G$,
then $G\in \mathfrak{F}$.
\end{Corollary}

\begin{Corollary}{\rm\cite[Corollary 3]{WH}}
Let $\mathfrak{F}$ be a saturated formation containing all supersoluble groups and $E$ a normal subgroup of $G$ with $G/E\in \mathfrak{F}$.
If all minimal subgroups
and all cyclic subgroups with order $4$ of all non-cyclic Sylow subgroups of $E$ are
$c$-normal in G,
then $G\in \mathfrak{F}$.
\end{Corollary}


\begin{thebibliography}{s2}


\bibitem{AN1}
Skiba, A.N.: On $\sigma$-subnormal and $\sigma$-permutable subgroups of finite groups. J. Algebra
 {\bf 436}, 1-16 (2015).

\bibitem{zs}
Skiba, A.N.: On some results in the theory of finite partially soluble groups.
Commun. Math. Stat. \textbf{4}, 281-309 (2016).


\bibitem{WS1}
Zhang, C., Guo, W., Liu, A-M.: On a generalisation of finite $T$-groups, Commun. Math. Stat. https://doi.org/10.1007/s40304-021-00240-z, (2021).


\bibitem{wei}
Wei, X.: On weakly $m$-$\s$-permutable subgroups of finite groups. Comm. Algebra \textbf{47}(3), 945-956 (2019).


\bibitem{Sch}
Schmidt, R.: \textit{Subgroup Lattices of Groups}. Walter de Gruyter, Berlin (1994).


\bibitem{H-permute}
Guo, W., Cao, C., Skiba, A.N., Sinitsa, D.A.: Finite groups with $\mathcal{H}$-permutable subgroups.
Commun. Math. Stat. {\bf 5}, 83-92 (2017).



\bibitem{AB2}
Ballester-Bolinches, A., Esteban-Romero, R., Asaad, M.: \textit{Products of Finite Groups}. Walter de Gruyter, Berlin, New York (2010).

\bibitem{W}
Guo, W.: \textit{Structure Theory for Canonical Classes of Finite Groups}. Springer, Heidelberg, New York, Dordrecht, London (2015).

\bibitem{MW}
Weinsten, M., et al.: \textit{Between Nilpotent and Soluble}. Polygonal Publishing House, Passaic (1982).


\bibitem{MSC}
Maier, R., Schmid, P.: The embedding of permutable subgroups in finite groups. Math. Z.
{\bf 131}, 269-272 (1973).

\bibitem{AN5}
Skiba, A. N.: On two questions of L.A. Shemetkov concerning hypercyclically embeded subgroups of finite groups. J. Group Theory {\bf 13}, 841-850 (2010).

\bibitem{BL}
Li, B.: On $\Pi$-property and $\Pi$-normality of subgroups of finite groups. J. Algebra {\bf 334}, 321-337 (2011).


\bibitem{WS2}
Guo, W., Skiba, A.N.: Finite groups with permutable complete Wielandt set of subgroups. J. Group Theory {\bf 18}, 191-200 (2015).


\bibitem{ZWG}
Zhang, C., Wu, Z., Guo, W.: On weakly $\s$-permutable subgroups of finite groups. Publ. Math. Debrecen \textbf{91}(3-4), 489-502 (2017).

\bibitem{MTH}
Hussain, M.T., Amjid, V.: Finite groups with weakly $m$-$\s$-permutable subgroups. J. Algebra Appl. {\bf 20}, 2150056 (15 pages) (2021).



\bibitem{WY1}
Wang, Y.: $C$-Normality of groups and its properties. J. Algebra
 {\bf 180}, 954-965 (1996).

\bibitem{Shaa}
Shaalan, A.: The influence of $\pi$-quasinormality of some subgroups on the structure of a finite group.
Acta Math. Hungar. {\bf 56}, 287-293 (1990).


\bibitem{SS}
Asaad, M.: On the solvability of finite groups. Arch. Math. {\bf 51}, 289-293 (1988).

\bibitem{Bu}
Buckley, J.: Finite groups whose minimal subgroups are normal. Math. Z. {\bf 116} 15-17 (1970).


\bibitem{Doerk}
Doerk, K., Hawkes, T.: \textit{Finite Soluble Groups}. Walter de Gruyter, Berlin (1992).

\bibitem{BNVS}
Knyagina, V.N., Monakhov, V.S.: On the $\pi'$-properties of a finite group possessing a Hall
$\pi$-subgroup. Sib. Math. J. \textbf{52}(2), 297-309 (2011).



\bibitem{GSO}
Guo, W., Skiba, A.N.: Finite groups with generalized Ore supplement conditions for
primary subgroups. J. Algebra {\bf 432}, 205-227 (2015).

\bibitem{Gfini}
Gorenstein, D.: \textit{Finite Groups}. Harper $\&$ Row Publishers, New York, Evanston, London (1968).

\bibitem{HU}
Huppert, B.: \textit{Endliche Gruppen \uppercase\expandafter{\romannumeral1}}. Springer-Verlag, Berlin (1967).

\bibitem{W1}
Guo, W.: \textit{The Theory of Classes of Groups}. Science Press-Kluwer Academic Publishers, Beijing, New York, Dordrecht, Boston, London (2000).


\bibitem{AN4}
Skiba, A. N.: A characterization of the hypercyclically embedded subgroups of finite groups. J. Pure Appl. Algebra {\bf 215}, 257-261 (2011).

\bibitem{AN3}
Skiba, A.N.: On weakly $s$-permutable subgroups of finite groups. J. Algebra {\bf 315}, 192-209 (2007).

\bibitem{HuNB}
Huppert, B., Blackburn, N.: \textit{Finite Groups. III}. Springer-Verlag, Berlin, New York (1982).

\bibitem{BaW}
Ballester-Bolinches, A., Wang, Y.: Finite groups with some $c$-normal minimal subgroups. J. Pure Appl. Algebra {\bf 153},
121-127 (2000).

\bibitem{BP}
 Ballester-Bolinches, A., Pedraza-Aguilera, M.C.: On minimal subgroups of finite groups. Acta Math. Hungar. {\bf 73},
335-342 (1996).

\bibitem{ACs}
Asaad, M., Cs$\ddot{o}$rg$\ddot{o}$, P.: The influence of minimal subgroups on the structure of finite group. Arch. Math. (Basel) {\bf 72}, 401-404 (1999).


\bibitem{LiWa}
Li, Y., Wang, Y.: The influence of minimal subgroups on the structure of a finite group. Proc. Amer. Math. Soc. {\bf 131},
337-341 (2002).

\bibitem{WH}
Wei, H.: On $c$-normal maximal and minimal subgroups of Sylow subgroups of
finite groups. Comm. Algebra {\bf 29}(5), 2193-2200 (2001).




\end{thebibliography}
\end{document}